\documentclass[11pt]{article}
\usepackage{epsfig}
\usepackage{amssymb}
\usepackage{amsmath,bm}
\usepackage{mathrsfs}
\newtheorem{theorem}{Theorem}[section]

\newtheorem{conjecture}[theorem]{Conjecture} 

\newcommand{\be}{\begin{equation}}
\newcommand{\ee}{\end{equation}}
\newcommand{\bea}{\begin{eqnarray}}
\newcommand{\eea}{\end{eqnarray}}

\usepackage{bm}
\usepackage{color}
\usepackage{graphicx}
\usepackage{cite}

\begin{document}


\title{A New Class of Relations for Homogeneous Symmetric Polynomials}
\author{Boris Y. Rubinstein\\
Stowers Institute for Medical Research
\\1000 50th St., Kansas City, MO 64110, U.S.A.}
\date{\today}

\maketitle
\begin{abstract}
Recently a new class of relations for the Bernoulli symmetric polynomials was introduced.
This manuscript shows that these relations are valid for arbitrary homogeneous symmetric polynomials.
\end{abstract}



\section{Introduction}
\label{relations_symm_poly}

Recently \cite{RubLinRel2025} we established existence of a new class of 
relations for the scalar partitions.
Their derivation is built on the
approach of J.J. Sylvester \cite[Lecture 5]{SylvLectures1859}
to the reduction of vector partition into the scalar ones and the 
corresponding theorem of A. Cayley for double partitions \cite{Cayley1860}.
Based on these findings we conjectured the following relation
\bea
&&\frac{B_{n}({\bf x}^m)}{\pi({\bf x}^m)} - \sum_{i=1}^m   
\frac{y_i^{m-n-1} B_{n}({\bf s}^m_i)}{\pi({\bf s}^m_i)} = 0, 
\quad
0 \le n \le m-1,
\label{Bernoulli_lin_polynomials} \\
&& s^m_{i,j} = y_i x_j - y_j x_i + x_i \delta_{ij}, \ (1\le j\le m),
\label{part_scalar_poly}
\eea
where $B_{n}({\bf x}^m)$ denote the homogeneous symmetric Bernoulli polynomials of higher order \cite{Bateman1953,Norlund1924},
${\bf x}^m=\{x_1,x_2,\ldots,x_m\}$
and $s^m_{i,j}$ are the elements of the rows ${\bf s}^m_i$ of the
$m \times m$ matrix ${\bf S}({\bf x}^m,{\bf y}^m) = {\bf y}^m \wedge {\bf x}^m + \mbox{diag}({\bf x}^m)$ 
given by the sum of the diagonal matrix with the elements ${\bf x}^m$ and the wedge product ${\bf y}^m \wedge {\bf x}^m$ of 
the vectors ${\bf x}^m$ and ${\bf y}^m$.
The direct computation shows that the relation (\ref{Bernoulli_lin_polynomials}) 
holds for $m \le 6$. 

In this manuscript we demonstrate that this relation is valid for the whole class of 
homogeneous symmetric polynomials that can be expressed by the complete Bell polynomials.
We then generalize it to {\it arbitrary} homogeneous symmetric polynomials.
In addition 
we 
describe the analog of the relation (\ref{Bernoulli_lin_polynomials}) 
for $n \ge m$.
We also report 
a new type of 
nonlinear relations for the Bernoulli numbers.

\section{Bernoulli polynomials as complete Bell polynomials}
\label{Bell}

Recall the generating functions $G_B(s,t)$ and $G_B(t,{\bf x}^m)$ for 
the univariate Bernoulli polynomial $B_n(s)$ and its 
multivariate symmetric counterpart -- the Bernoulli polynomial of higher order \cite{Norlund1924,Rub04}
\be
G_B(s,t) = \frac{t e^{st}}{e^{t}-1} = \sum_{n=0}^{\infty} B_n(s)\frac{t^n}{n!},
\quad
G_B(t,{\bf x}^m) = \prod_{i=1}^m \frac{x_i t}{e^{x_it}-1} =
\sum_{n=0}^{\infty} B_n({\bf x}^m)
\frac{t^{n}}{n!}\;.
\label{BernoulliGF}
\ee
Use $G_B(t,{\bf x}^m)$ together with 
\be
\ln G_B(0,x t) = \ln \frac{xt}{e^{xt}-1} =
\sum_{k=1}^{\infty} (-1)^{k-1} \frac{x^k B_k}{k} \frac{t^k}{k!},
\label{BernoulliGFLog}
\ee
to obtain 
$$
\ln G_B(t,{\bf x}^m) = \sum_{k=1}^{\infty} (-1)^{k-1} \frac{p_k({\bf  x}^m) B_k}{k} \frac{t^k}{k!},
\quad
\quad
p_k({\bf  x}^m) = \sum_{i=1}^m x_i^k,
$$
where we introduce the power sums $p_k({\bf  x}^m)$.
This expression establishes the relation between the Bernoulli polynomials of higher order
and the complete Bell polynomials \cite{Rub2009}
\be
B_n({\bf  x}^m) = \mathbb B_n(b_1,b_2,\ldots),
\quad
b_k = (-1)^{k-1}B_k p_k({\bf  x}^m)/k,
\label{BernoulliBell}
\ee
where 
the complete Bell polynomials $ \mathbb B_n(b_1,b_2,\ldots)$ are defined by the generating function \cite{Riordan}
\be
\exp \left(\sum_{i=1}^{\infty}\frac{b_i}{i!}t^i \right) =
\sum_{n=0}^{\infty} \frac{\mathbb B_n(b_1,b_2,\ldots)}{n!}t^n.
\label{BellGF}
\ee
The expression for $ \mathbb B_{n}({\bf b})\equiv \mathbb B_{n}(b_1,b_2,\ldots)$ depends
on the first $n$ elements of the vector ${\bf b}=\{b_1,b_2,\ldots\}$ only.
Write down the explicit expressions for $ \mathbb B_{n}({\bf b})$ for small $n\le 5$:
\bea
&& \mathbb B_{0}({\bf b}) = 1,
\quad\quad
 \mathbb B_{1}({\bf b}) = b_1,
\quad\quad
\mathbb B_{2}({\bf b}) =  b_1^2+ b_2,
\quad
\mathbb B_{3}({\bf b}) =  b_1^3+ 3b_1b_2,
\nonumber\\
&&  \mathbb B_{4}({\bf b}) =  b_1^4+ 6b_1^2b_2+ 3b_2^2+b_4, 
\quad
\mathbb B_{5}({\bf b}) =  b_1^5+ 10b_1^3b_2 + 15b_1b_2^2 + 5b_1b_4,
 \label{Bell_explicit} 
\eea
where we take into account that $B_{2k+1}=b_{2k+1}=0$ for $k > 0$.

Use (\ref{BernoulliBell}) to rewrite relations (\ref{Bernoulli_lin_polynomials}) through the 
complete Bell polynomials 
\bea
&&\frac{\mathbb B_{n}(b_1,b_2,\ldots)}{\pi({\bf x}^m)} - \sum_{i=1}^m   
\frac{y_i^{m-n-1} \mathbb B_{n}(b_{1i},b_{2i},\ldots)}{\pi({\bf s}^m_i)} = 0, 
\label{linrel_BernBell_polynomials} \\
&&b_k = (-1)^{k-1}B_k p_k({\bf  x}^m)/k,
\quad
b_{ki} = (-1)^{k-1}B_k p_k({\bf s}^m_i)/k.
\nonumber
\eea

\section{Relations for complete Bell symmetric polynomials}
\label{BellSymm}

It is instructive to generalize the relation (\ref{linrel_BernBell_polynomials}) to the 
case of symmetric polynomials expressed through the complete Bell polynomials. 
Consider a family of symmetric polynomials defined as follows.
Start with a polynomial $F_n(s)$ of a single variable $s$  
defined by an exponential generating function $G(s,t)$ and 
similar to (\ref{BernoulliGF},\ref{BernoulliGFLog})
generalize it 
to a symmetric  polynomial $F_n({\bf x}^m)$ of $m$ variables ${\bf x}^m$ through 
the multiplicative generating function $G(t,{\bf x}^m)$ where
\be
G(s,t) = \sum_{n=0}^{\infty} F_n(s)\frac{t^n}{n!},
\quad
G(t,{\bf x}^m) =
\prod_{i=1}^m G(s_0, x_i t) = \sum_{n=0}^{\infty}F_n({\bf x}^m) \frac{t^n}{n!}.
\label{GF_symm_poly_general}
\ee
Let the general form of the Taylor expansion of $\ln G(s_0,t x)$ at the fixed value $s=s_0$ read
\be
\ln G(s_0,t x)  = \sum_{k=1}^{\infty} \frac{t^k}{k!} a_{k} x^{k},
\label{ln_poly_general}
\ee
then we find for $G(t,{\bf x}^m)$ 
\be
\ln G(t,{\bf x}^m)  = \sum_{k=1}^{\infty} a_{k} p_{k}({\bf x}^m) \frac{t^k}{k!}.
\label{ln_symm_poly_general}
\ee
This relation leads to the corresponding complete Bell polynomials 
\be
G(t,{\bf x}^m)  = \exp 
\left(
\sum_{k=1}^{\infty} a_{k} p_{k}({\bf x}^m)\frac{t^k}{k!} 
\right) = 
\sum_{n=0}^{\infty} \mathbb B_n(f_1,f_2,\ldots)\frac{t^n}{n!},
\label{symm_poly_general_log}
\ee
and comparing it to the second equation in (\ref{GF_symm_poly_general}) we 
obtain
\be
F_n({\bf x}^m) = \mathbb B_n(f_1,f_2,\ldots),
\quad
f_k = a_{k} p_{k}({\bf x}^m). 
\label{symm_poly_general}
\ee
Write down the explicit expressions for $F_n({\bf x}^m)$ for small $n\le 5$:
\bea
&& F_0({\bf x}^m) = 1,
\quad\quad
F_1({\bf x}^m) = f_1,
\quad\quad
F_2({\bf x}^m) =  f_1^2+ f_2,
\nonumber\\
&&
F_3({\bf x}^m) =  f_1^3+ 3f_1f_2 +f_3,
\quad  
F_4({\bf x}^m) =  f_1^4+ 6f_1^2f_2+ 3f_2^2 + 4f_1f_3+f_4,
\nonumber\\
&&
F_5({\bf x}^m) =  f_1^5+ 10f_1^3f_2 + 15f_1f_2^2 + 10f_1^2f_3 + 10f_2f_3 
+ 5f_1f_4 + f_5.
 \label{Gen_poly_explicit} 
\eea

\begin{conjecture}
\label{conj1}
Introduce a function
\be
U_n({\bf x}^m,{\bf y}^m) = V_n({\bf x}^m) - 
\sum_{i=1}^m y_i^{m-n-1}  V_n({\bf s}^m_i),
\quad
 V_n({\bf x}^m) = F_{n}({\bf x}^m)/\pi({\bf x}^m),
\label{main_G_def}
\ee
where the polynomials $F_{n}({\bf x}^m)=\alpha_n \mathbb B_n(f_1,f_2,\ldots)$ 
with the factor $\alpha_n$ independent of ${\bf x}^m$ and $f_k$
defined in (\ref{symm_poly_general}); ${\bf s}^m_i$ are the rows of the matrix ${\bf S}({\bf x}^m,{\bf y}^m)$
with the elements given in (\ref{part_scalar_poly}).
The following relation holds
\be
U_n({\bf x}^m,{\bf y}^m) = 0,
\quad
0 \le n \le m-1.
\label{linrel_General_polynomials} 
\ee
\end{conjecture}

It is attractive to expand the range of $n$ values in (\ref{linrel_General_polynomials}) beyond 
$n=m-1$ expecting to arrive at $U_n({\bf x}^m,{\bf y}^m)$ being 
a rational function in ${\bf x}^m,{\bf y}^m$.   
Direct computation shows however that it has a simpler structure
$$U_{n}({\bf x}^m,{\bf y}^m) = R_{n-m}({\bf x}^m,{\bf y}^m)/\pi^{n-m+1}({\bf y}^m),
\quad n \ge m,$$
where $R_{n-m}({\bf x}^m,{\bf y}^m)$ is a polynomial 
of the degree $(n-m)$ in ${\bf x}^m$ and degree $m(n-m+1)$ in ${\bf y}^m$.
Setting ${\bf y}^m={\bf 1}^m$ we find
that for the Bernoulli polynomials of higher order 
$R_{n-m}({\bf x}^m,{\bf 1}^m)=0$, while 
for the other polynomials including $T_n({\bf x}^m)$ introduced in \cite{Fel21} and 
thoroughly investigated in \cite{FelOno25} it is a symmetric polynomial in ${\bf x}^m$.
These results lead to 
\begin{conjecture}
\label{conj2}
For $n \ge m$ the function $U_n({\bf x}^m,{\bf y}^m)$ defined 
in (\ref{main_G_def}) at ${\bf y}^m={\bf 1}^m$ is a homogeneous 
symmetric polynomial in ${\bf x}^m$ of the degree  $n-m$
\be
U_n({\bf x}^m,{\bf 1}^m)= R_{n-m}({\bf x}^m,{\bf 1}^m) = Z_{n-m}({\bf x}^m),
\quad
n \ge m.
\label{linrel_General_polynomials_nonzero} 
\ee
\end{conjecture} 
The examples of the homogeneous symmetric polynomials $F_{n}({\bf x}^m)=\alpha_n \mathbb B_n(f_1,f_2,\ldots)$ 
satisfying (\ref{linrel_General_polynomials},\ref{linrel_General_polynomials_nonzero})
are given in the table below.
\vskip0.2cm
\renewcommand{\arraystretch}{1.5}
\begin{tabular}{|l|c|c|c|c|}  
\hline 
Polynomial & $G(s,t)$ & $a_k$ & $\alpha_n$  & $s_0$ \\
\hline
Legendre $P_n$  &  $e^{s t}J_0(t \sqrt{1-s^2})$ & 
$a_{2k+1} = 0, \ a_{2k} = \frac{d^{2k}\ln J_0(t)}{dt^{2k}}\vert_{t=0}$  & 1 & 0\\ 
\hline
Laguerre $L_n$ &  $e^{-ts/(1-t)}/(1-t)$ & 
$a_{k} = (k-1)!$   & $1/n!$  & 0\\
\hline
Hermite $H_n$ & $e^{2 s t - t^2}$ & 
$a_2 = -2,\ a_i=0,\ i\ne 2$ & 1 & 0\\
\hline
Fibonacci $F_n$ & $1/(1-xt-t^2)$ &
$a_{2k+1}=0,\ a_{2k} = 2((2k-1)!)$ & $1/n!$  & 0\\
\hline
Bernoulli $B_n$ & $t e^{s t}/(e^t-1)$ &
$a_k = (-1)^{k-1} B_k/k$   &   1&  0\\
\hline
$T_n$ & $(e^t-1)/(t e^{s t})$ &
$a_k = (-1)^k B_k/k$  &  1 & 0 \\
\hline
Euler $E_n$ & $2 e^{s t}/(e^t+1)$ &
$a_1=-1/2, a_k = E_{k-1}(0)/2$  &  1 & 0 \\
\hline
Bell ${\cal B}_n$ & $e^{(e^t-1)s}$   &
$a_k = 1$   & 1 & 1 \\
\hline
\end{tabular}

\section{Structure of polynomials $Z_{n}({\bf x}^m)$}
\label{Z explicit}

From the definition of $Z_{n}({\bf x}^m)$
its expansion 
follows
\be
Z_{n}({\bf x}^m) = \sum_{\bf k} z^m_{\bf k} P_{n,\bf k}({\bf x}^m),
\quad
P_{n,\bf k}({\bf x}^m)=
p_1^{k_1}({\bf x}^m) 
p_2^{k_2}({\bf x}^m)
\ldots
p_n^{k_n}({\bf x}^m),
\quad
\sum_{i=1}^n i k_i = n.
\label{Z_poly_expansion}
\ee
where ${\bf k} = \{k_1,k_2,\ldots,k_n\}$ denotes the integer vector of 
powers $k_i$ of the power sums $p_i({\bf x}^m)$.
Start from $Z_{0}({\bf x}^m)=z_0^m$ and $Z_{1}({\bf x}^m) = z^m_{\{1\}} p_1({\bf x}^m)$
with
\begin{eqnarray}
&&z^2_0=a_1^2+3a_2, \
z^3_0=2(4a_1^3+12a_1a_2+a_3), \
z^4_0=81a_1^4+270a_1^2a_2+75a_2^2+36a_1a_3+5a_4, 
\nonumber \\
&&z^2_{\{1\}}=2(a_1^3+3a_1a_2+a_3), \
z^3_{\{1\}}=16a_1(a_1^3+3a_1a_2+a_3), \
\label{Z0_Z1} \\
&&z^4_{\{1\}}=2(81a_1^5+270a_1^3a_2+75a_1a_2^2+90a_1^2a_3+30a_2a_3+5a_1a_4).
\nonumber
\end{eqnarray}
Continue with $n=2$ where $Z_{2}({\bf x}^m) = z^m_{\{2,0\}} p_1^2({\bf x}^m) + z^m_{\{0,1\}} p_2({\bf x}^m)$
having
\begin{eqnarray}
&&
z^2_{\{2,0\}}=(5a_1^4+6a_1^2a_2-21a_2^2-4a_1a_3-5a_4)/2,
\
z^2_{\{0,1\}}=(a_1^4+30a_1^2a_2+63a_2^2+28a_1a_3+15a_4)/2, 
\nonumber \\
&&
z^3_{\{2,0\}}=16a_1^5-120a_1a_2^2-20a_1^2a_3-60a_2a_3-20a_1a_4-4a_5,
\nonumber \\
&&z^3_{\{0,1\}}=16a_1^5+240a_1^3a_2+480a_1a_2^2+160a_1^2a_3+180a_2a_3+80a_1a_4+11a_5, \
\nonumber \\
&&z^4_{\{2,0\}}=(243a_1^6-405a_1^4a_2
-3285a_1^2a_2^2-825a_2^3-540a_1^3a_3-1980a_1a_2a_3
\nonumber \\
&&
-90a_3^2-435a_1^2a_4-285a_2a_4-102a_1a_5-7a_6)/2,
\label{Z2} \\
&&
z^4_{\{0,1\}}=(729a_1^6+8505a_1^4a_2
+18225a_1^2a_2^2+4125a_2^3+4860a_1^3a_3+9180a_1a_2a_3
\nonumber \\
&&
+450a_3^2+2295a_1^2a_4+1425a_2a_4+414a_1a_5+35a_6)/2.
\nonumber 
\end{eqnarray}
The expressions for the coefficients $z^m_{\bf k}$ in $Z_{n}({\bf x}^m)$ for $n=3,4$ 
are quite cumbersome and presented in the Appendix \ref{Z_4}. 

\section{Nonlinear relations for the Bernoulli numbers}
\label{Bern number relation}

It is instructive to obtain the conditions on $a_i$ values in (\ref{symm_poly_general}) for which
the polynomials $Z_n({\bf x}^m)$ vanish. Use the expressions for $z^m_{\bf k}$ 
in (\ref{Z0_Z1},\ref{Z2}) as well as in the Appendix \ref{Z_4} and set them to zero. Solving these equations to find 
$a_i,\ i>1$ expressed through $a_1$ we obtain
\be
a_{3}=a_{5}=a_{7}=\ldots = 0,
\quad
a_2=-\frac{a_1^2}{3},
\quad
a_4=\frac{2a_1^4}{15},
\quad
a_6=-\frac{16a_1^6}{63},
\quad
a_8=\frac{16a_1^8}{15},
\quad
\ldots
\label{a_values}
\ee 
We observe that the values of
$a_k$ in (\ref{a_values}) satisfy $a_k = - (2 a_1)^{k} B_k/k$ leading to
the corresponding symmetric polynomials 
$F_n({\bf x}^m) =  (-2a_1)^n B_n({\bf x}^m)$ and for these polynomials the 
relations $Z_n({\bf x}^m)=0$ hold.
This means that the Bernoulli numbers $B_n$ satisfy a large (possible infinite) set of {\it nonlinear} relations
examples of which are presented below
\begin{eqnarray}
&&
B_1^2 - 3 B_2/2 = 0,
\quad
B_1^4 - 15 B_1^2 B_2 + 63 B_2^2/4-15 B_4/4 = 0,
\nonumber \\
&&
B_1^4 + B_1^2 B_2 - 9B_2^2/4 + 5B_4/4 = 0,
\nonumber \\
&&
10 B_1^6 - 135 B_1^4 B_2 - 135 B_1^2 B_2^2/2 + 585 B_2^3/4 - 
 75 B_1^2 B_4/2 - 405 B_2 B_4/4 + 7 B_6 = 0,
\label{BernNumberRels} \\
&&
 8 B_1^6 + 21 B_1^4 B_2 + 105 B_1^2 B_2^2 - 105 B_2^3/4 + 
 35 B_1^2 B_4 + 315 B_2 B_4/4 - 28 B_6/3 = 0.
 \nonumber
\end{eqnarray}
It should be underlined that the relations shown above are not unique. 
For example, instead of the last two relations
we might use also
$$
8 B_1^6 - 294 B_1^4 B_2 + 2100 B_1^2 B_2^2 - 
 2205 B_2^3/2 - 560 B_1^2 B_4 + 2835 B_2 B_4/2 -  140 B_6 = 0.
$$
Write down the explicit expressions for the symmetric polynomial $B_n({\bf x}^m)$
for $n \le 6$ (to save space we drop the argument ${\bf x}^m$ of the power sums):
\bea
&& B_0 = 1, \quad
B_1 = -p_1/2, \quad
B_2 = (3p_1^2-p_2)/12,\quad
B_3 =  -p_1(p_1^2-p_2)/8, 
\label{Bern_explicit} \\
&& 
B_4 =(15 p_1^4 - 30 p_1^2 p_2 + 5 p_2^2 + 2 p_4)/240,
\quad
B_5 = -p_1 (3 p_1^4 - 10 p_1^2 p_2 + 5 p_2^2 + 2 p_4)/96,
\nonumber \\
&& 
B_6 = (63 p_1^6 - 315 p_1^4 p_2 + 315 p_1^2 p_2^2 - 35 p_2^3 + 
 126 p_1^2 p_4 - 42 p_2 p_4 - 16 p_6)/4032.
\nonumber
\eea
\section{Relation for power sum products}
\label{relations_power_sums}

Assuming that the conjectures \ref{conj1} and \ref{conj2} are valid 
the definition of the symmetric polynomial (\ref{symm_poly_general}) 
implies more fundamental relations for the power sums products.
To show that note that $F_{n}({\bf x}^m)$ equals the 
sum of the specific power sums or their products with {\it arbitrary} 
coefficients similar to (\ref{Z_poly_expansion}), namely 
\be
F_{n}({\bf x}^m) = \sum_{\bf k} A_{\bf k}({\bf a}) P_{n,\bf k}({\bf x}^m),
\quad
P_{n,\bf k}({\bf x}^m)=
p_1^{k_1}({\bf x}^m) 
p_2^{k_2}({\bf x}^m)
\ldots
p_n^{k_n}({\bf x}^m),
\quad
\sum_{i=1}^n i k_i = n.
\label{G_poly_expansion}
\ee
Substitute this expansion into (\ref{symm_poly_general}) and observe
that similar relation must hold for the function 
\be
U_{n,\bf k}({\bf x}^m,{\bf y}^m) = V_{n,\bf k}({\bf x}^m) - 
\sum_{i=1}^m y_i^{m-n-1}  V_{n,\bf k}({\bf s}^m_i),
\quad
 V_{n,\bf k}({\bf x}^m) = P_{n,\bf k}({\bf x}^m)/\pi({\bf x}^m),
\label{main_powersum_def}
\ee
Similarly for $n \ge m$ we find 
$U_{n,\bf k}({\bf x}^m,{\bf 1}^m) = Y_{n-m,\bf k}({\bf x}^m)$.

\begin{conjecture}
\label{conj1powersum}
For an arbitrary homogeneous symmetric polynomial $S_n({\bf x}^m)$ of the degree $n$
expressed through the power sums $p_l({\bf x}^m)$
\be
S_{n}({\bf x}^m) = \sum_{\bf k} C_{n,\bf k} P_{n,\bf k}({\bf x}^m),
\quad
P_{n,\bf k}({\bf x}^m)=
p_1^{k_1}({\bf x}^m) 
p_2^{k_2}({\bf x}^m)
\ldots
p_n^{k_n}({\bf x}^m),
\quad
\sum_{i=1}^n i k_i = n,
\label{S_poly_expansion}
\ee
define function
\be
U_{n}({\bf x}^m,{\bf y}^m) = V_{n}({\bf x}^m) - 
\sum_{i=1}^m y_i^{m-n-1}  V_{n}({\bf s}^m_i),
\quad
 V_{n}({\bf x}^m) = S_{n}({\bf x}^m)/\pi({\bf x}^m),
\label{main_S_def}
\ee
where ${\bf s}^m_i$ has the elements $s^m_{i,j} = y_i x_j - y_j x_i + x_i \delta_{ij}$.
The following relations hold
\bea
&& U_{n}({\bf x}^m,{\bf y}^m) = 0,
\quad
0 \le n \le m-1, 
\label{linrel_S_polynomials} \\
&& 
U_{n}({\bf x}^m,{\bf y}^m) =R_{n-m}({\bf x}^m,{\bf y}^m)/\pi^{n-m+1}({\bf y}^m),
\quad
 n \ge m,
 \label{linrel_S_polynomials_nonzero0}
\eea
where 
$R_{n-m}({\bf x}^m,{\bf y}^m)$ is the homogeneous polynomial of the degree $(n-m)$ in ${\bf x}^m$ 
with the coefficients being polynomials of the degree $m(n-m+1)$ in ${\bf y}^m$;
the polynomials $R_{n-m}({\bf x}^m,{\bf y}^m)$
satisfy an additional symmetry, namely, 
they are invariant under any permutations of the columns of the following $2 \times m$ matrix
$$
\left(
\begin{array}{ccccc}
x_1 & x_2 & x_3 & \ldots & x_m \\
y_1 & y_2 & y_3 & \ldots & y_m
\end{array}
\right).
$$
Setting in (\ref{linrel_S_polynomials_nonzero0}) ${\bf y}^m={\bf 1}^m$ we obtain
\be
U_{n}({\bf x}^m,{\bf 1}^m) = R_{n-m}({\bf x}^m,{\bf 1}^m)= Y_{n-m}({\bf x}^m),
\quad
n \ge m,
\label{linrel_S_polynomials_nonzero} 
\ee
where $Y_{n-m}({\bf x}^m)$ denotes the homogeneous symmetric polynomial of the degree $(n-m)$ 
computed using the expansion
into the 
polynomials $Y_{n-m,\bf k}({\bf x}^m)$
\be
Y_{n-m}({\bf x}^m) = \sum_{\bf k} C_{n,\bf k} Y_{n-m,\bf k}({\bf x}^m).
\label{Y_basis} 
\ee
\end{conjecture} 
We find that $Y_{n-1,\bf k}({\bf x}^1)=0$ for all $n$, and 
the explicit expressions of $Y_{n-m,\bf k}({\bf x}^m)$ for $2 \le n \le 8$ are
presented in Appendices \ref{Ybasis_m=2}, \ref{Ybasis_m=3}, \ref{Ybasis_m=4}
with $m=2,3,4,$ respectively. 

In Section \ref{Bern number relation} we note that among all homogeneous symmetric polynomials 
of the degree $n$ expressible through the complete Bell polynomials only the Bernoulli polynomials 
force all $Y_{n-m}({\bf x}^m)$ to vanish. 
It is worth to find out if there exist values of the coefficients $C_{n,\bf k}$
for which polynomials $Y_{n-m}({\bf x}^m)$ defined in (\ref{Y_basis}) turn to zero 
for fixed $n$ and all $2 \le m \le n$.

To answer this question we first have to find the number of the coefficients $C_{n,\bf k}$
to be determined -- it is equal to the number of integer vectors ${\bf k}$ of the length $n$ 
with the nonnegative elements $k_i$ satisfying the condition $\sum_{i=1}^n i k_i = n$
and it is just the unrestricted partition $P(n)$ of the integer $n$.

The products $P_{n-m,\bf k}({\bf x}^m)$ in the equation (\ref{Y_basis}) for the polynomial $Y_{n-m}({\bf x}^m)$ 
contain only the power sums $p_j({\bf x}^m)$ with $1 \le j \le m$ and their number is given by
the partition $P_m(n-m)$ of $(n-m)$ into parts less or equal $m$
that satisfies the recurrence relation $P_m(n) = P_{m-1}(n)+P_m(n-m)$.
Using this relation we find the number $E_n$ of the linear equations for $C_{n,\bf k}$
given by the sum
$$
E_n = \sum_{m=2}^n P_m(n-m)= P_n(n) - P_1(n) = P(n) - 1.
$$
Thus for each value of $n$ there exists a family of coefficients $C_{n,\bf k}$ defining the polynomials $\bar S_{n}({\bf x}^m)$ for which 
the relations $Y_{n-m}({\bf x}^m)=0$ hold for all $2 \le m \le n$.

The solutions for $C_{n,\bf k}$ for $n \le 6$ are shown in Appendix \ref{coeff_C} and using them in 
(\ref{S_poly_expansion}) we obtain the corresponding polynomials $\bar S_{n}({\bf x}^m)$ in the form
\bea
\bar S_{1}({\bf x}^m) &=& C_{1, \{1\}} p_1,\quad 
\bar S_{2}({\bf x}^m) = C_{2, \{2, 0\}} (p_1^2 - p_2/3),
\quad
\bar S_{3}({\bf x}^m) = C_{3, \{3, 0, 0\}} p_1 (p_1^2 - p_2),
\nonumber\\
 \bar S_{4}({\bf x}^m) &=& C_{4, \{4, 0, 0, 0\}} (p_1^4  - 3p_2^2 - 4  p_1 p_3 + 6  p_4)
 + C_{4, \{2, 1, 0, 0\}} (p_1^2 p_2 -5p_2^2/3 -2 p_1 p_3  + 44p_4/15),
\nonumber  \\
\bar S_{5}({\bf x}^m) &=&  C_{5, \{5, 0, 0, 0, 0\}}( p_1^5 - 25 p_1 p_2^2/3  - 10 p_1^2 p_3+ 50 p_2 p_3/3+ 62 p_1 p_4/3 - 20  p_5)
\label{S_n_1to5} \\
 &+& C_{5, \{3, 1, 0, 0, 0\}}( p_1^3 p_2  -3 p_1 p_2^2  -3 p_1^2 p_3 +  5 p_2 p_3 + 6  p_1 p_4 -6  p_5),
\quad
p_i \equiv p_i({\bf x}^m).
\nonumber
\eea
Setting here the following values for $C_{n,\bf k}$
\begin{eqnarray*}
&&
C_{1, \{1\}} = -1/2,\ 
C_{2, \{2, 0\}} = 1/4, \
C_{3, \{3, 0, 0\}} = -1/8, \ 
 C_{4, \{4, 0, 0, 0\}} = 1/16, \\
 && 
C_{4, \{2, 1, 0, 0\}} = -1/8,\ 
C_{5, \{5, 0, 0, 0, 0\}}= -1/32, \ 
C_{5, \{3, 1, 0, 0, 0\}} = 5/48,
\end{eqnarray*}
we recover as a particular case the symmetric Bernoulli polynomials (\ref{Bern_explicit}).



\newpage
{\LARGE \bf Appendices}

\appendix

\section{Expressions for $Z_{n}({\bf x}^m)$ with $n=3,4$}
\renewcommand{\theequation}{\thesection\arabic{equation}}
\setcounter{equation}{0}
\label{Z_4}

We have $Z_{3}({\bf x}^m) = z^m_{\{3,0,0\}} p_1^3({\bf x}^m) +  z^m_{\{1,1,0\}} p_1({\bf x}^m)p_2({\bf x}^m) +
z^m_{\{0,0,1\}} p_3({\bf x}^m)$,
where
\begin{eqnarray}
&&
z^2_{\{3,0,0\}}=5(a_1^5-2a_1^3a_2-9a_1a_2^2-8a_1^2a_3-8a_2a_3-5a_1a_4-a_5)/2,
\nonumber \\
&&z^2_{\{1,1,0\}}=(3a_1^5+90a_1^3a_2+165a_1a_2^2+120a_1^2a_3+120a_2a_3+65a_1a_4+13a_5)/2, 
\ z^2_{\{0,0,1\}} = 0,
\nonumber \\
&&
z^3_{\{3,0,0\}}=
(16a_1^6-480a_1^4a_2
-1440a_1^2a_2^2-120a_2^3-400a_1^3a_3-840a_1a_2a_3
\nonumber \\
&&
-110a_3^2-120a_1^2a_4+60a_2a_4+6a_1a_5+7a_6)/3,
\nonumber \\
&&z^3_{\{1,1,0\}}=
2(32a_1^6+480a_1^4a_2
+1080a_1^2a_2^2+120a_2^3+280a_1^3a_3+480a_1a_2a_3
\nonumber \\
&&
+50a_3^2+60a_1^2a_4-60a_2a_4-18a_1a_5-7a_6),
\label{Z3} \\
&&z^3_{\{0,0,1\}}=
(-64a_1^6-960a_1^4a_2
-2880a_1^2a_2^2-9605a_2^3+160a_1^3a_3+480a_1a_2a_3
\nonumber \\
&&
+80a_3^2+480a_1^2a_4+480a_2a_4+336a_1a_5+56a_6)/3,
\nonumber \\
&&z^4_{\{3,0,0\}}=-81a_1^7 -2457 a_1^5 a_2 - 6825 a_1^3 a_2^2 - 2205 a_1a_2^3 -
 1575 a_1^4 a_3 - 4410 a_1^2 a_2 a_3 -  525 a_2^2 a_3
 \nonumber \\
 && - 630 a_1 a_3^2 - 455 a_1^3 a_4 - 105 a_1 a_2 a_4 - 
 35 a_3 a_4 - 21 a_1^2 a_5 + 63 a_2 a_5 + 21 a_1 a_6 + 3 a_7,
 \nonumber  \\
&&
z^4_{\{1,1,0\}}=1458 a_1^7 + 17010 a_1^5 a_2 + 39690 a_1^3 a_2^2 + 
 13650 a_1 a_2^3 + 8505 a_1^4 a_3 + 20790 a_1^2 a_2 a_3 + 
 2835 a_2^2 a_3 
\nonumber \\ 
 &&+ 2520 a_1 a_3^2 + 1890 a_1^3 a_4 - 
 210 a_1 a_2 a_4 + 105 a_3 a_4 - 252 a_1^2 a_5 - 
 420 a_2 a_5 - 182 a_1 a_6 - 19 a_7,
\nonumber \\
 &&
z^4_{\{0,0,1\}}=-729 a_1^7 - 8505 a_1^5 a_2 - 23625 a_1^3 a_2^2 - 
 13125 a_1 a_2^3 + 630 a_1 a_3^2 + 2205 a_1^3 a_4 + 
 3675 a_1 a_2 a_4 
\nonumber \\
 &&+ 420 a_3 a_4 + 1701 a_1^2 a_5 + 
 945 a_2 a_5 + 385 a_1 a_6 + 34 a_7.
\nonumber 
\end{eqnarray}

For $n=4$ we have 
\bea
Z_{4}({\bf x}^m) &=& z^m_{\{4,0,0,0\}} p_1^4({\bf x}^m) 
+ z^m_{\{2,1,0,0\}} p_1^2({\bf x}^m)p_2({\bf x}^m) +z^m_{\{0,2,0,0\}} p_2^2({\bf x}^m)
\nonumber \\
&+&
z^m_{\{1,0,1,0\}} p_1({\bf x}^m)p_3({\bf x}^m)+z^m_{\{0,0,0,1\}} p_4({\bf x}^m),
\nonumber
\eea
where $z^2_{\{1,0,1,0\}}=z^2_{\{0,0,0,1\}}=z^3_{\{0,0,0,1\}}=0$.
We have for $m=2$ three nonzero coefficients
\begin{eqnarray}
&&
z^2_{\{4,0,0,0\}}=(9 a_1^6 - 45 a_1^4 a_2 + 45 a_1^2 a_2^2 + 195 a_2^3 - 
 180 a_1^3 a_3 + 180 a_1 a_2 a_3 + 90 a_3^2
 \nonumber \\&&- 105 a_1^2 a_4 + 
 135 a_2 a_4 - 6 a_1 a_5 + 7 a_6)/4,
\label{n=4,m=2} \\
&&z^2_{\{2,1,0,0\}}=(10 a_1^6 + 270 a_1^4 a_2 - 270 a_1^2 a_2^2 - 1170 a_2^3 + 
 440 a_1^3 a_3 - 1080 a_1 a_2 a_3 - 380 a_3^2 
  \nonumber \\&&+ 
 150 a_1^2 a_4 - 810 a_2 a_4 - 60 a_1 a_5 - 42 a_6)/4, 
 \nonumber \\&&
 z^2_{\{0,2,0,0\}} = (a_1^6 + 75 a_1^4 a_2 + 1125 a_1^2 a_2^2 + 1755 a_2^3 + 
 140 a_1^3 a_3 + 2100 a_1 a_2 a_3 + 490 a_3^2  
 \nonumber \\&&+ 
 255 a_1^2 a_4 + 1215 a_2 a_4 + 186 a_1 a_5 + 63 a_6)/4.
\nonumber 
\end{eqnarray}
For $m=3$ we find
\begin{eqnarray}
&&
z^3_{\{4,0,0,0\}}=(-48 a_1^7 - 1680 a_1^5 a_2 - 1680 a_1^3 a_2^2 + 
 5040 a_1 a_2^3 - 1260 a_1^4 a_3  
+ 840 a_1^2 a_2 a_3 + 3360 a_2^2 a_3  
 \nonumber \\&&+ 840 a_1 a_3^2 + 560 a_1^3 a_4 + 
 3360 a_1 a_2 a_4 + 770 a_3 a_4 + 672 a_1^2 a_5 + 
 630 a_2 a_5 + 252 a_1 a_6 + 29 a_7)/6,
\label{n=4,m=3} \\
&&z^3_{\{2,1,0,0\}}=(672 a_1^7 + 8064 a_1^5 a_2 - 35280 a_1 a_2^3 + 
 3360 a_1^4 a_3 - 17640 a_1^2 a_2 a_3 - 22680 a_2^2 a_3 - 
 7140 a_1 a_3^2
  \nonumber \\&& - 5040 a_1^3 a_4 - 22680 a_1 a_2 a_4 - 
 5040 a_3 a_4 - 4158 a_1^2 a_5 - 3906 a_2 a_5 - 
 1470 a_1 a_6 - 168 a_7)/6, 
 \nonumber \\&&
 z^3_{\{0,2,0,0\}} = (192 a_1^7 + 6048 a_1^5 a_2 + 45360 a_1^3 a_2^2 + 
 65520 a_1 a_2^3 + 6720 a_1^4 a_3 + 60480 a_1^2 a_2 a_3 + 
 41580 a_2^2 a_3
  \nonumber \\&& + 13440 a_1 a_3^2 + 7560 a_1^3 a_4 + 
 32760 a_1 a_2 a_4 + 6930 a_3 a_4 + 4032 a_1^2 a_5 + 
 5418 a_2 a_5 + 1344 a_1 a_6 + 171 a_7)/6,
\nonumber \\
&&
 z^3_{\{1,0,1,0\}} = (-384 a_1^7 - 5376 a_1^5 a_2 - 13440 a_1^3 a_2^2 - 
 3360 a_1^2 a_2 a_3 - 3360 a_2^2 a_3 + 1680 a_1 a_3^2 + 
 4480 a_1^3 a_4
 \nonumber \\&& + 6720 a_1 a_2 a_4 + 1120 a_3 a_4 + 
 2856 a_1^2 a_5 + 504 a_2 a_5 + 840 a_1 a_6 + 64 a_7)/6.
\nonumber 
\end{eqnarray}
Finally for $m=4$ we obtain
\begin{eqnarray}
&&
z^4_{\{4,0,0,0\}}=5(-2997 a_1^8 - 41580 a_1^6 a_2 - 83370 a_1^4 a_2^2 + 
 18060 a_1^2 a_2^3 + 12915 a_2^4 - 20664 a_1^5 a_3
\nonumber \\&& - 18480 a_1^3 a_2 a_3 + 59640 a_1 a_2^2 a_3 + 
 7560 a_1^2 a_3^2 + 14280 a_2 a_3^2 + 490 a_1^4 a_4 + 
 26460 a_1^2 a_2 a_4 + 9870 a_2^2 a_4
 \nonumber \\&& + 12040 a_1 a_3 a_4 + 
 665 a_4^2 + 2408 a_1^3 a_5 + 5880 a_1 a_2 a_5 + 
 1736 a_3 a_5 + 532 a_1^2 a_6 + 84 a_2 a_6 - 8 a_1 a_7 - 
 9 a_8)/8,
 \nonumber \\
&&z^4_{\{2,1,0,0\}}=(24786 a_1^8 + 258552 a_1^6 a_2 + 396900 a_1^4 a_2^2 - 
 321720 a_1^2 a_2^3 - 129150 a_2^4 + 81648 a_1^5 a_3 
\nonumber \\&&- 
 151200 a_1^3 a_2 a_3 - 589680 a_1 a_2^2 a_3 - 
 115920 a_1^2 a_3^2 - 115920 a_2 a_3^2 - 34020 a_1^4 a_4 - 
 271320 a_1^2 a_2 a_4
 \nonumber \\&& - 98700 a_2^2 a_4 - 
 109200 a_1 a_3 a_4 - 6650 a_4^2 - 23184 a_1^3 a_5 - 
 48048 a_1 a_2 a_5 - 12432 a_3 a_5 - 3976 a_1^2 a_6 
\nonumber \\&&- 
 840 a_2 a_6 + 272 a_1 a_7 + 90 a_8)/8, 
\label{n=4,m=4} \\&&
 z^4_{\{0,2,0,0\}} = (6561 a_1^8 + 183708 a_1^6 a_2 + 1241730 a_1^4 a_2^2 + 
 2060100 a_1^2 a_2^3 + 401625 a_2^4 + 204120 a_1^5 a_3
 \nonumber \\&& + 
 1678320 a_1^3 a_2 a_3 + 1912680 a_1 a_2^2 a_3 + 
 385560 a_1^2 a_3^2 + 264600 a_2 a_3^2 + 164430 a_1^4 a_4 + 
 805140 a_1^2 a_2 a_4
 \nonumber \\&& + 303450 a_2^2 a_4 + 
 244440 a_1 a_3 a_4 + 19355 a_4^2 + 65016 a_1^3 a_5 + 
 154728 a_1 a_2 a_5 + 22680 a_3 a_5 + 16380 a_1^2 a_6
 \nonumber \\&& + 
 13020 a_2 a_6 + 2088 a_1 a_7 + 117 a_8)/8,
\nonumber \\
&&
 z^4_{\{1,0,1,0\}} = (-29160 a_1^8 - 381024 a_1^6 a_2 - 1285200 a_1^4 a_2^2 - 
 1092000 a_1^2 a_2^3 - 105000 a_2^4 - 108864 a_1^5 a_3 
 \nonumber \\&&- 
 604800 a_1^3 a_2 a_3 - 504000 a_1 a_2^2 a_3 - 
 20160 a_1^2 a_3^2 - 20160 a_2 a_3^2 - 25200 a_1^4 a_4 - 
 178080 a_1^2 a_2 a_4
 \nonumber \\&& - 75600 a_2^2 a_4 - 6720 a_1 a_3 a_4 - 
 3640 a_4^2 - 20160 a_1^3 a_5 - 81984 a_1 a_2 a_5 - 
 1344 a_3 a_5 - 16352 a_1^2 a_6
 \nonumber \\&& - 14560 a_2 a_6 - 
 4800 a_1 a_7 - 456 a_8)/8,
\nonumber \\
&&
 z^4_{\{0,0,0,1\}} = (13122 a_1^8 + 204120 a_1^6 a_2 + 850500 a_1^4 a_2^2 + 
 945000 a_1^2 a_2^3 + 131250 a_2^4 + 81648 a_1^5 a_3
 \nonumber \\&& + 
 453600 a_1^3 a_2 a_3 + 378000 a_1 a_2^2 a_3 + 
 45360 a_1^2 a_3^2 + 25200 a_2 a_3^2 + 102060 a_1^4 a_4 + 
 340200 a_1^2 a_2 a_4
 \nonumber \\&& + 94500 a_2^2 a_4 + 45360 a_1 a_3 a_4 + 
 4550 a_4^2 + 69552 a_1^3 a_5 + 115920 a_1 a_2 a_5 + 
 7728 a_3 a_5 + 32760 a_1^2 a_6
 \nonumber \\&& + 18200 a_2 a_6 + 
 6864 a_1 a_7 + 570 a_8)/8.
 \nonumber 
\end{eqnarray}

\section{Expressions for $Y_{n-2,\bf k}({\bf x}^2)$}
\renewcommand{\theequation}{\thesection\arabic{equation}}
\setcounter{equation}{0}
\label{Ybasis_m=2}
We present the expressions for  $Y_{n-m,\bf k}=Y_{n-m,\bf k}({\bf x}^m)$ for $m=2$ and $2 \le n \le 8$.
To save space we use $p_k$  instead of $p_k({\bf x}^2)$.
\begin{eqnarray}
&&Y_{0,\{2,0\}}=1,\
Y_{0,\{0,1\}}=3,\
Y_{1,\{3,0,0\}}=Y_{1,\{1,1,0\}}=Y_{1,\{0,0,1\}}=2p_1,
\nonumber \\&&
Y_{2,\{4,0,0,0\}}=\frac{5p_1^2+p_2}{2},\
Y_{2,\{2,1,0,0\}}=\frac{p_1^2+5p_2}{2},\
Y_{2,\{0,2,0,0\}}=\frac{-7p_1^2+21p_2}{2},
\nonumber \\&&
Y_{2,\{1,0,1,0\}}=\frac{-p_1^2+7p_2}{2},\
Y_{2,\{0,0,0,1\}}=\frac{-5p_1^2+3p_2}{2}, 
\nonumber \\&&
Y_{3,\{5, 0, 0, 0, 0\}}=\frac{p_1(5p_1^2 + 3p_2)}{2},\
Y_{3,\{3, 1, 0, 0, 0\}}=\frac{p_1(-p_1^2 + 9p_2)}{2},\
Y_{3,\{1, 2, 0, 0, 0\}}=\frac{p_1(-3p_1^2 + 11p_2)}{2},
\nonumber \\&&
Y_{3,\{2, 0, 1, 0, 0\}}=Y_{3,\{0, 1, 1, 0, 0\}}=2p_1(-p_1^2 + 3p_2),\
Y_{3,\{1, 0, 0, 1, 0\}}=
Y_{3,\{0, 0, 0, 0, 1\}}=\frac{p_1(-5p_1^2 +13p_2)}{2},
\nonumber
\end{eqnarray}
\begin{eqnarray}
&&
Y_{4,\{6, 0, 0, 0, 0, 0\}}=\frac{(p_1^2 + p_2)(9p_1^2 + p_2)}{4},\
Y_{4,\{4, 1, 0, 0, 0, 0\}}=\frac{-3p_1^4 + 18p_1^2pp_2 + 5p_2^2}{4},\
\nonumber \\&&
Y_{4,\{2, 2, 0, 0, 0, 0\}}=\frac{p_1^4 - 6p_1^2p_2 + 25p_2^2}{4},\
Y_{4,\{0, 3, 0, 0, 0, 0\}}=\frac{13(p_1^2 - 3p_2)^2}{4},\
\nonumber \\&&
Y_{4,\{3, 0, 1, 0, 0, 0\}}=\frac{-9p_1^4 + 22p_1^2p_2 + 7p_2^2}{4},\
Y_{4,\{1, 1, 1, 0, 0, 0\}}=\frac{3p_1^4 - 18p_1^2p_2 + 35p_2^2}{4},\
\nonumber \\&&
Y_{4,\{0, 0, 2, 0, 0, 0\}}=\frac{9p_1^4 - 38p_1^2p_2 + 49p_2^2}{4},\
Y_{4,\{2, 0, 0, 1, 0, 0\}}=-\frac{(7p_1^2 - 17p_2)(p_1^2 + p_2)}{4},\
\nonumber \\&&
Y_{4,\{0, 1, 0, 1, 0, 0\}}=\frac{9(p_1^2 - 3p_2)^2}{4},\
Y_{4,\{1, 0, 0, 0, 1, 0\}}=\frac{-p_1^4 - 10p_1^2p_2 + 31p_2^2}{4},\
\nonumber \\&&
Y_{4,\{0, 0, 0, 0, 0, 1\}}=\frac{7(p_1^2 - 3p_2)^2}{4},
\nonumber
\end{eqnarray}
\begin{eqnarray}
&&
Y_{5,\{7, 0, 0, 0, 0, 0, 0\}}=p_1(p_1^2 + p_2)(2p_1^2 + p_2),\
Y_{5,\{5, 1, 0, 0, 0, 0, 0\}}=-\frac{p_1(p_1^2 - 7p_2)(p_1^2 + p_2)}{2},
\nonumber \\&&
Y_{5,\{3, 2, 0, 0, 0, 0, 0\}}=p_1(p_1^4 - 5p_1^2p_2 + 10p_2^2),\
Y_{5,\{1, 3, 0, 0, 0, 0, 0\}}=\frac{p_1(p_1^2 - 7p_2)(p_1^2 - 3p_2)}{2},
\nonumber \\&&
Y_{5,\{4, 0, 1, 0, 0, 0, 0\}}=-\frac{p_1(7p_1^2 - 19p_2)(p_1^2 + pp_2)}{4},\
Y_{5,\{2, 1, 1, 0, 0, 0, 0\}}=\frac{p_1(7p_1^4 - 36p_1^2p_2 + 53p_2^2)}{4},
\nonumber \\&&
Y_{5,\{0, 2, 1, 0, 0, 0, 0\}}=\frac{p_1(p_1^4 - 20p_1^2p_2 + 43p_2^2)}{4},\
Y_{5,\{1, 0, 2, 0, 0, 0, 0\}}=\frac{p_1(7p_1^4 - 30p_1^2p_2 + 35p_2^2)}{2},
\nonumber \\&&
Y_{5,\{3, 0, 0, 1, 0, 0, 0\}}=-p_1(p_1^4 + p_1^2p_2 - 8p_2^2),\
Y_{5,\{1, 1, 0, 1, 0, 0, 0\}}=\frac{3p_1(p_1^2 - 3p_2)^2}{2},
\nonumber \\&&
Y_{5,\{0, 0, 1, 1, 0, 0, 0\}}=\frac{p_1(11p_1^4 - 52p_1^2p_2 + 65p_2^2)}{4},
\nonumber \\&&
Y_{5,\{2, 0, 0, 0, 1, 0, 0\}}=Z_{5,\{0, 1, 0, 0, 1, 0, 0\}}=\frac{p_1(p_1^2 - 7p_2)(3p_1^2 - 7p_2)}{4},
\nonumber \\&&
Y_{5,\{1, 0, 0, 0, 0, 1, 0\}}=Y_{5,\{0, 0, 0, 0, 0, 0, 1\}}=p_1(2p_1^2 - 5p_2)(p_1^2 - 3p_2),
\nonumber
\end{eqnarray}
\begin{eqnarray}
&&
Y_{6,\{8, 0, 0, 0, 0, 0, 0, 0\}}=\frac{15p_1^6 + 23p_1^4p_2 + 17p_1^2p_2^2 + p_2^3}{8},\
Y_{6,\{6, 1, 0, 0, 0, 0, 0, 0\}}=\frac{-p_1^6 + 11p_1^4p_2 + 41p_1^2p_2^2 + 5p_2^3}{8},\
\nonumber \\&&
Y_{6,\{4, 2, 0, 0, 0, 0, 0, 0\}}=\frac{7p_1^6 - 33p_1^4p_2 + 57p_1^2p_2^2 + 25p_2^3}{8},\
Y_{6,\{2, 3, 0, 0, 0, 0, 0, 0\}}=\frac{-9p_1^6 + 51p_1^4p_2 - 111p_1^2p_2^2 + 125p_2^3}{8},\
\nonumber \\&&
Y_{6,\{0, 4, 0, 0, 0, 0, 0, 0\}}=\frac{-17p_1^6 + 183p_1^4p_2 - 591p_1^2p_2^2 + 609p_2^3}{8},\
Y_{6,\{5, 0, 1, 0, 0, 0, 0, 0\}}=\frac{-9p_1^6 + 5p_1^4p_2 + 53p_1^2p_2^2 + 7p_2^3}{8},\
\nonumber \\&&
Y_{6,\{3, 1, 1, 0, 0, 0, 0, 0\}}=\frac{11p_1^6 - 55p_1^4p_2 + 65p_1^2p_2^2 + 35p_2^3}{8},\
Y_{6,\{1, 2, 1, 0, 0, 0, 0, 0\}}=\frac{-17p_1^6 + 93p_1^4p_2 - 195p_1^2p_2^2 + 175p_2^3}{8},\
\nonumber \\&&
Y_{6,\{2, 0, 2, 0, 0, 0, 0, 0\}}=\frac{21p_1^6 - 85p_1^4p_2 + 71p_1^2p_2^2 + 49p_2^3}{8},\
Y_{6,\{0, 1, 2, 0, 0, 0, 0, 0\}}=\frac{-31p_1^6 + 167p_1^4p_2 - 325p_1^2p_2^2 + 245p_2^3}{8},\
\nonumber \\&&
Y_{6,\{4, 0, 0, 1, 0, 0, 0, 0\}}=\frac{-5p_1^6 - 17p_1^4p_2 + 61p_1^2p_2^2 + 17p_2^3}{8},\
Y_{6,\{2, 1, 0, 1, 0, 0, 0, 0\}}=\frac{3p_1^6 - 13pp_1^4p_2 - 19p_1^2p_2^2 + 85p_2^3}{8},\
\nonumber \\&&
Y_{6,\{0, 2, 0, 1, 0, 0, 0, 0\}}=-\frac{3(7p_1^6 - 53p_1^4p_2 + 145p_1^2p_2^2 - 139p_2^3)}{8},\
Y_{6,\{1, 0, 1, 1, 0, 0, 0, 0\}}=\frac{7p_1^6 - 11p_1^4p_2 - 59p_1^2p_2^2 + 119p_2^3}{8},\
\nonumber \\&&
Y_{6,\{0, 0, 0, 2, 0, 0, 0, 0\}}=\frac{-5p_1^6 + 75p_1^4p_2 - 267p_1^2p_2^2 + 285p_2^3}{8},\
Y_{6,\{3, 0, 0, 0, 1, 0, 0, 0\}}=\frac{5p_1^6 - 47p_1^4p_2 + 67p_1^2p_2^2 + 31p_2^3}{8},\
\nonumber \\&&
Y_{6,\{1, 1, 0, 0, 1, 0, 0, 0\}}=\frac{-11p_1^6 + 61p_1^4p_2 - 149p_1^2p_2^2 + 155p_2^3}{8},\
 \nonumber \\&&
Y_{6,\{0, 0, 1, 0, 1, 0, 0, 0\}}=\frac{-19p_1^6 + 115p_1^4p_2 - 257p_1^2p_2^2 + 217p_2^3}{8},\
\nonumber \\&&
Y_{6,\{2, 0, 0, 0, 0, 1, 0, 0\}}=\frac{9p_1^6 - 45p_1^4p_2 + 27p_1^2p_2^2 + 65p_2^3}{8},\
Y_{6,\{0, 1, 0, 0, 0, 1, 0, 0\}}=\frac{-23p_1^6 + 147p_1^4p_2 - 357p_1^2p_2^2 + 321p_2^3}{8},\
\nonumber \\&&
Y_{6,\{1, 0, 0, 0, 0, 0, 1, 0\}}=\frac{p_1^6 + 9p_1^4p_2 - 81p_1^2p_2^2 + 127p_2^3}{8},\
Y_{6,\{0, 0, 0, 0, 0, 0, 0, 1\}}=-\frac{3(5p_1^6 - 35p_1^4p_2 + 91p_1^2p_2^2 - 85p_2^3)}{8}.
\nonumber 
\end{eqnarray}

\section{Expressions for $Y_{n-3,\bf k}({\bf x}^3)$}
\renewcommand{\theequation}{\thesection\arabic{equation}}
\setcounter{equation}{0}
\label{Ybasis_m=3}
The explicit expressions for $Y_{n-m,\bf k}=Y_{n-m,\bf k}({\bf x}^m)$ for $m=3$ and $3 \le n \le 8$ are shown below.
To save space we use $p_k$  instead of $p_k({\bf x}^3)$.
\begin{eqnarray}
&&
Y_{0,\{3,0,0\}}=Y_{0,\{1,1,0\}}=8, Y_{0,\{0,0,1\}}=2,\
\nonumber \\&&
Y_{1,\{4,0,0,0\}}=16p_1,\   
Y_{1,\{2,1,0,0\}}=8p_1,\    
Y_{1,\{1,0,1,0\}}=4p_1,  
Y_{1,\{0,2,0,0\}}=Y_{1,\{0,0,0,1\}}=0, 
\nonumber \\&&
Y_{2,\{5, 0, 0, 0, 0\}}=16 (p_1^2 + p_2),\
Y_{2,\{3, 1, 0, 0, 0\}}=24 p_2,\
Y_{2,\{1, 2, 0, 0, 0\}}=-8 (p_1^2 - 4 p_2),
\nonumber \\&&
Y_{2,\{2, 0, 1, 0, 0\}}=-2 (p_1^2 - 8 p_2),\
Y_{2,\{0, 1, 1, 0, 0\}}=-6 (p_1^2 - 3 p_2),
\nonumber \\&&
Y_{2,\{1, 0, 0, 1, 0\}}=-4 (p_1^2 - 4 p_2),\
Y_{2,\{0, 0, 0, 0, 1\}}=-4 p_1^2 + 11 p_2,
\nonumber
\end{eqnarray}
\begin{eqnarray}
&&
Y_{3,\{6, 0, 0, 0, 0, 0\}}=\frac{16 (p_1^3 + 12 p_1 p_2 - 4 p_3)}{3},\
Y_{3,\{4, 1, 0, 0, 0, 0\}}=-\frac{32 (p_1^3 - 6 p_1 p_2 + 2 p_3)}{3},\
\nonumber \\&&
Y_{3,\{2, 2, 0, 0, 0, 0\}}=-\frac{16 (2 p_1^3 - 9 p_1 p_2 + 4 p_3)}{3},\
Y_{3,\{0, 3, 0, 0, 0, 0\}}=-\frac{8 (p_1^3 - 6 p_1 p_2 + 8 p_3)}{3},\
\nonumber \\&&
Y_{3,\{3, 0, 1, 0, 0, 0\}}=-\frac{4 (5 p_1^3 - 21 p_1 p_2 - 2 p_3)}{3},\
Y_{3,\{1, 1, 1, 0, 0, 0\}}=-\frac{2 (7 p_1^3 - 24 p_1 p_2 - 4 p_3)}{3},\
\nonumber \\&&
Y_{3,\{0, 0, 2, 0, 0, 0\}}=-\frac{11 p_1^3 + 30 p_1 p_2 + 8 p_3}{3},\
Y_{3,\{2, 0, 0, 1, 0, 0\}}=-\frac{8 (p_1^3 - 3 p_1 p_2 - 4 p_3)}{3},\
\nonumber \\&&
Y_{3,\{0, 1, 0, 1, 0, 0\}}=\frac{4 (p_1^3 - 6 p_1 p_2 + 8 p_3)}{3},\
Y_{3,\{1, 0, 0, 0, 1, 0\}}=\frac{p_1^3 - 18 p_1 p_2 + 56 p_3}{3},\
\nonumber \\&&
Y_{3,\{0, 0, 0, 0, 0, 1\}}=\frac{7 (p_1^3 - 6 p_1 p_2 + 8 p_3)}{3},
\nonumber
\end{eqnarray}
\begin{eqnarray}
&&
Y_{4,\{7, 0, 0, 0, 0, 0, 0\}}=-8 (p_1^4 - 14 p_1^2 p_2 - 4 p_2^2 + 8 p_1 p_3),\
Y_{4,\{5, 1, 0, 0, 0, 0, 0\}}=-\frac{8 (5 p_1^4 - 24 p_1^2 p_2 - 18 p_2^2 + 16 p_1 p_3)}{3},
\nonumber \\&&
Y_{4,\{3, 2, 0, 0, 0, 0, 0\}}=-\frac{8 (p_1^4 - 27 p_2^2 + 8 p_1 p_3)}{3},\
Y_{4,\{1, 3, 0, 0, 0, 0, 0\}}=8 (p_1^4 - 7 p_1^2 p_2 + 13 p_2^2),
\nonumber \\&&
Y_{4,\{4, 0, 1, 0, 0, 0, 0\}}=-2 (p_1^2 - 4 p_2)(3 p_1^2 + 4 p_2),\
Y_{4,\{2, 1, 1, 0, 0, 0, 0\}}=\frac{2 (p_1^4 - 21 p_1^2 p_2 + 72 p_2^2 - 4 p_1 p_3)}{3},
\nonumber \\&&
Y_{4,\{0, 2, 1, 0, 0, 0, 0\}}=\frac{2 (8 p_1^4 - 54 p_1^2 p_2 + 99 p_2^2 - 8 p_1 p_3)}{3},\
Y_{4,\{1, 0, 2, 0, 0, 0, 0\}}=2 p_1^4 - 17 p_1^2 p_2 + 32 p_2^2 + 4 p_1 p_3,
\nonumber \\&&
Y_{4,\{3, 0, 0, 1, 0, 0, 0\}}=\frac{4 (2 p_1^4 - 18 p_1^2 p_2 + 27 p_2^2 + 16 p_1 p_3)}{3},\
Y_{4,\{1, 1, 0, 1, 0, 0, 0\}}=\frac{4 (4 p_1^4 - 27 p_1^2 p_2 + 39 p_2^2 + 8 p_1 p_3)}{3},
\nonumber \\&&
Y_{4,\{0, 0, 1, 1, 0, 0, 0\}}=\frac{11 p_1^4 - 72 p_1^2 p_2 + 99 p_2^2 + 16 p_1 p_3}{3},\
Y_{4,\{2, 0, 0, 0, 1, 0, 0\}}=\frac{16 p_1^4 - 99 p_1^2 p_2 + 96 p_2^2 + 68 p_1 p_3}{3},
\nonumber \\&&
Y_{4,\{0, 1, 0, 0, 1, 0, 0\}}=5 p_1^4 - 31 p_1^2 p_2 + 43 p_2^2 + 4 p_1 p_3,\
Y_{4,\{1, 0, 0, 0, 0, 1, 0\}}=6 p_1^4 - 35 p_1^2 p_2 + 32 p_2^2 + 20 p_1 p_3,
\nonumber \\&&
Y_{4,\{0, 0, 0, 0, 0, 0, 1\}}=\frac{29 p_1^4 - 168 p_1^2 p_2 + 171 p_2^2 + 64 p_1 p_3}{6},
\nonumber
\end{eqnarray}
\begin{eqnarray}
&&
Y_{5,\{8, 0, 0, 0, 0, 0, 0, 0\}}=-\frac{32 (p_1^5 - 8 p_1^3 p_2 - 18 p_1 p_2^2 + 8 p_1^2 p_3 + 8 p_2 p_3}{3},\
\nonumber \\&&
Y_{5,\{6, 1, 0, 0, 0, 0, 0, 0\}}=-\frac{8 (p_1^5 + 8 p_1^3 p_2 - 84 p_1 p_2^2 + 8 p_1^2 p_3 + 40 p_2 p_3}{3},\
\nonumber \\&&
Y_{5,\{4, 2, 0, 0, 0, 0, 0, 0\}}=\frac{16 (2 p_1^5 - 18 p_1^3 p_2 + 45 p_1 p_2^2 + 4 p_1^2 p_3 - 24 p_2 p_3)}{3},\
\nonumber \\&&
Y_{5,\{2, 3, 0, 0, 0, 0, 0, 0\}}=\frac{8 (5 p_1^5 - 40 p_1^3 p_2 + 81 p_1 p_2^2 + 16 p_1^2 p_3 - 56 p_2 p_3)}{3},\
\nonumber \\&&
Y_{5,\{0, 4, 0, 0, 0, 0, 0, 0\}}=\frac{16 (p_1^2 - 4 p_2) (p_1^3 - 6 p_1 p_2 + 8 p_3)}{3},\
\nonumber \\&&
Y_{5,\{5, 0, 1, 0, 0, 0, 0, 0\}}=\frac{4 (p_1^5 - 26 p_1^3 p_2 + 90 p_1 p_2^2 + 8 p_1^2 p_3 - 28 p_2 p_3)}{3},\
\nonumber \\&&
Y_{5,\{3, 1, 1, 0, 0, 0, 0, 0\}}=\frac{2 (11 p_1^5 - 86 p_1^3 p_2 + 174 p_1 p_2^2 + 16 p_1^2 p_3 - 52 p_2 p_3)}{3},\
\nonumber \\&&
Y_{5,\{1, 2, 1, 0, 0, 0, 0, 0\}}=\frac{4 (5 p_1^5 - 36 p_1^3 p_2 + 63 p_1 p_2^2 + 10 p_1^2 p_3 -  24 p_2 p_3)}{3},\
\nonumber \\&&
Y_{5,\{2, 0, 2, 0, 0, 0, 0, 0\}}=\frac{13 p_1^5 - 86 p_1^3 p_2 + 144 p_1 p_2^2 -  4 p_1^2 p_3 + 32 p_2 p_3}{3},\
\nonumber \\&&
Y_{5,\{0, 1, 2, 0, 0, 0, 0, 0\}}=\frac{7 p_1^5 - 46 p_1^3 p_2 +  78 p_1 p_2^2 - 16 p_1^2 p_3 + 40 p_2 p_3}{3},\
\nonumber \\&&
Y_{5,\{4, 0, 0, 1, 0, 0, 0, 0\}}=\frac{8 (3 p_1^5 - 22 p_1^3 p_2 + 33 p_1 p_2^2 + 12 p_1^2 p_3 -  8 p_2 p_3)}{3},\
\nonumber \\&&
Y_{5,\{2, 1, 0, 1, 0, 0, 0, 0\}}=\frac{4 (4 p_1^5 - 28 p_1^3 p_2 + 45 p_1 p_2^2 +  8 p_1^2 p_3 - 8 p_2 p_3)}{3},\
Y_{5,\{0, 2, 0, 1, 0, 0, 0, 0\}}=0,\
\nonumber \\&&
Y_{5,\{1, 0, 1, 1, 0, 0, 0, 0\}}=\frac{2 (3 p_1^5 - 16 p_1^3 p_2 + 15 p_1 p_2^2 - 6 p_1^2 p_3 +  40 p_2 p_3)}{3},\
\nonumber \\&&
Y_{5,\{0, 0, 0, 2, 0, 0, 0, 0\}}=-\frac{4 (p_1^2 - 4 p_2) (p_1^3 - 6 p_1 p_2 + 8 p_3)}{3},\
\nonumber \\&&
Y_{5,\{3, 0, 0, 0, 1, 0, 0, 0\}}=\frac{23 p_1^5 - 134 p_1^3 p_2 + 126 p_1 p_2^2 +  64 p_1^2 p_3 + 44 p_2 p_3}{3},\
\nonumber \\&&
Y_{5,\{1, 1, 0, 0, 1, 0, 0, 0\}}=\frac{7 p_1^5 - 34 p_1^3 p_2 + 12 p_1 p_2^2 - 4 p_1^2 p_3 +  100 p_2 p_3}{3},\
\nonumber \\&&
Y_{5,\{0, 0, 1, 0, 1, 0, 0, 0\}}=\frac{-p_1^5 + 16 p_1^3 p_2 - 45 p_1 p_2^2 - 38 p_1^2 p_3 + 128 p_2 p_3}{3},\
\nonumber \\&&
Y_{5,\{2, 0, 0, 0, 0, 1, 0, 0\}}=\frac{13 p_1^5 - 62 p_1^3 p_2 + 18 p_1 p_2^2 + 20 p_1^2 p_3 +  104 p_2 p_3}{3},\
\nonumber \\&&
Y_{5,\{0, 1, 0, 0, 0, 1, 0, 0\}}=-\frac{5 (p_1^2 - 4 p_2) (p_1^3 - 6 p_1 p_2 +  8 p_3)}{3},\
\nonumber \\&&
Y_{5,\{1, 0, 0, 0, 0, 0, 1, 0\}}=\frac{3 p_1^5 + 4 p_1^3 p_2 - 75 p_1 p_2^2 - 18 p_1^2 p_3 + 164 p_2 p_3}{3},\
\nonumber \\&&
Y_{5,\{0, 0, 0, 0, 0, 0, 0, 1\}}=-2 (p_1^2 - 4 p_2) (p_1^3 - 6 p_1 p_2 + 8 p_3).
\nonumber 
\end{eqnarray}

\section{Expressions for $Y_{n-4,\bf k}({\bf x}^4)$}
\renewcommand{\theequation}{\thesection\arabic{equation}}
\setcounter{equation}{0}
\label{Ybasis_m=4}
Below shown the expressions for $Y_{n-m,\bf k}=Y_{n-m,\bf k}({\bf x}^m)$ for $m=4$ and $4 \le n \le 8$.
To save space we use $p_k$  instead of $p_k({\bf x}^4)$.
\begin{eqnarray}
&&
Y_{0,\{4,0,0,0\}}=81,\
Y_{0,\{2,1,0,0\}}=45,\
Y_{0,\{0,2,0,0\}}=25,
Y_{0,\{1,0,1,0\}}=9,\
Y_{0,\{0,0,0,1\}}=5, 
\nonumber \\&&
Y_{1,\{5, 0, 0, 0, 0\}}=162 p_1,\
Y_{1,\{3, 1, 0, 0, 0\}}= 54 p_1,\
Y_{1,\{1, 2, 0, 0, 0\}}= 10 p_1,\
Y_{1,\{2, 0, 1, 0, 0\}}=18 p_1,
\nonumber \\&&
Y_{3,\{0, 1, 1, 0, 0\}}=6 p_1,\
Y_{1,\{1, 0, 0, 1, 0\}}=2 p_1,\
Y_{1,\{0, 0, 0, 0, 1\}}=2 p_1,
\nonumber
\end{eqnarray}
\begin{eqnarray}
&&
Y_{2,\{6, 0, 0, 0, 0, 0\}}=\frac{243 (p_1^2 + 3 p_2)}{2},\
Y_{2,\{4, 1, 0, 0, 0, 0\}}=-\frac{27 (p_1^2 - 21 p_2)}{2},\
\nonumber \\&&
Y_{2,\{2, 2, 0, 0, 0, 0\}}=-\frac{73 p_1^2 + 405 p_2}{2},\
Y_{2,\{0, 3, 0, 0, 0, 0\}}=-\frac{55 (p_1^2 - 5 p_2)}{2},\
\nonumber \\&&
Y_{2,\{3, 0, 1, 0, 0, 0\}}=-\frac{27 (p_1^2 - 9 p_2)}{2},\
Y_{2,\{1, 1, 1, 0, 0, 0\}}=-\frac{3 (11 p_1^2 - 51 p_2)}{2},\
\nonumber \\&&
Y_{2,\{0, 0, 2, 0, 0, 0\}}=-\frac{9 (p_1^2 - 5 p_2)}{2},\
Y_{2,\{2, 0, 0, 1, 0, 0\}}=\frac{-29 p_1^2 + 153 p_2}{2},\
\nonumber \\&&
Y_{2,\{0, 1, 0, 1, 0, 0\}}=-\frac{19 (p_1^2 - 5 p_2)}{2},\
Y_{2,\{1, 0, 0, 0, 1, 0\}}=-\frac{17 p_1^2 + 69 p_2}{2},\
\nonumber \\&&
Y_{2,\{0, 0, 0, 0, 0, 1\}}=-\frac{7 (p_1^2 - 5 p_2)}{2},
\nonumber
\end{eqnarray}
\begin{eqnarray}
&&
Y_{3,\{7, 0, 0, 0, 0, 0, 0\}}=-81 (p_1^3 - 18 p_1 p_2 + 9 p_3),\
Y_{3,\{5, 1, 0, 0, 0, 0, 0\}}=-9 (13 p_1^3 - 90 p_1 p_2 + 45 p_3),
\nonumber \\&&
Y_{3,\{3, 2, 0, 0, 0, 0, 0\}}=-65 p_1^3 + 378 p_1 p_2 - 225 p_3,\
Y_{3,\{1, 3, 0, 0, 0, 0, 0\}}=-21 p_1^3 + 130 p_1 p_2 - 125 p_3,
\nonumber \\&&
Y_{3,\{4, 0, 1, 0, 0, 0, 0\}}=-9 p_1 (5 p_1^2 - 27 p_2),\
Y_{3,\{2, 1, 1, 0, 0, 0, 0\}}=-3 p_1 (7 p_1^2 - 33 p_2),
\nonumber \\&&
Y_{3,\{0, 2, 1, 0, 0, 0, 0\}}=-p_1 (5 p_1^2 - 27 p_2),\
Y_{3,\{1, 0, 2, 0, 0, 0, 0\}}=-9 (p_1^3 - 4 p_1 p_2 - p_3),
\nonumber \\&&
Y_{3,\{3, 0, 0, 1, 0, 0, 0\}}=-13 p_1^3 + 54 p_1 p_2 + 63 p_3,\
Y_{3,\{1, 1, 0, 1, 0, 0, 0\}}= -p_1^3 - 2 p_1 p_2 + 35 p_3,
\nonumber \\&&
Y_{3,\{0, 0, 1, 1, 0, 0, 0\}}=-p_1^3 + 3 p_1 p_2 + 12 p_3,\
Y_{3,\{2, 0, 0, 0, 1, 0, 0\}}= -p_1^3 - 12 p_1 p_2 + 81 p_3,
\nonumber \\&&
Y_{3,\{0, 1, 0, 0, 1, 0, 0\}}=3 p_1^3 - 20 p_1 p_2 + 45 p_3,\
Y_{3,\{1, 0, 0, 0, 0, 1, 0\}}=3 p_1^3 - 26 p_1 p_2 + 55 p_3,
\nonumber \\&&
Y_{3,\{0, 0, 0, 0, 0, 0, 1\}}=3 p_1^3 - 19 p_1 p_2 + 34 p_3,
\nonumber
\end{eqnarray}
\begin{eqnarray}
&&
Y_{4,\{8, 0, 0, 0, 0, 0, 0, 0\}}=-\frac{81 (37 p_1^4 - 306 p_1^2 p_2 - 81 p_2^2 + 360 p_1 p_3 - 162 p_4)}{8},\
\nonumber \\&&
Y_{4,\{6, 1, 0, 0, 0, 0, 0, 0\}}=-\frac{27 (55 p_1^4 - 342 p_1^2 p_2 - 243 p_2^2 + 504 p_1 p_3 - 270 p_4)}{8},\
\nonumber \\&&
Y_{4,\{4, 2, 0, 0, 0, 0, 0, 0\}}=\frac{-397 p_1^4 + 1890 p_1^2 p_2 + 5913 p_2^2 - 6120 p_1 p_3 + 4050 p_4}{8},\
\nonumber \\&&
Y_{4,\{2, 3, 0, 0, 0, 0, 0, 0\}}=\frac{43 p_1^4 - 766 p_1^2 p_2 + 4905 p_2^2 - 2600 p_1 p_3 + 2250 p_4}{8},\
\nonumber \\&&
Y_{4,\{0, 4, 0, 0, 0, 0, 0, 0\}}=\frac{123 p_1^4 - 1230 p_1^2 p_2 + 3825 p_2^2 - 1000 p_1 p_3 + 1250 p_4}{8},\
\nonumber \\&&
Y_{4,\{5, 0, 1, 0, 0, 0, 0, 0\}}=-\frac{9 (41 p_1^4 - 162 p_1^2 p_2 - 405 p_2^2 + 216 p_1 p_3 - 162 p_4)}{8},\
\nonumber \\&&
Y_{4,\{3, 1, 1, 0, 0, 0, 0, 0\}}=-\frac{3 (11 p_1^4 + 90 p_1^2 p_2 - 999 p_2^2 + 360 p_1 p_3 - 270 p_4)}{8},\
\nonumber \\&&
Y_{4,\{1, 2, 1, 0, 0, 0, 0, 0\}}=\frac{71 p_1^4 - 702 p_1^2 p_2 + 2277 p_2^2 -  600 p_1 p_3 + 450 p_4}{8},\
\nonumber \\&&
Y_{4,\{2, 0, 2, 0, 0, 0, 0, 0\}}=\frac{9 (3 p_1^4 - 46 p_1^2 p_2 + 153 p_2^2 - 8 p_1 p_3 + 18 p_4)}{8},\
\nonumber \\&&
Y_{4,\{0, 1, 2, 0, 0, 0, 0, 0\}}=\frac{3 (17 p_1^4 - 138 p_1^2 p_2 + 315 p_2^2 - 24 p_1 p_3 + 30 p_4)}{8},\
\nonumber \\&&
Y_{4,\{4, 0, 0, 1, 0, 0, 0, 0\}}=\frac{7 p_1^4 - 486 p_1^2 p_2 + 2349 p_2^2 - 360 p_1 p_3 +  1458 p_4}{8},\
\nonumber \\&&
Y_{4,\{2, 1, 0, 1, 0, 0, 0, 0\}}=\frac{63 p_1^4 - 646 p_1^2 p_2 + 1917 p_2^2 -  424 p_1 p_3 + 810 p_4}{8},\
\nonumber \\&&
Y_{4,\{0, 2, 0, 1, 0, 0, 0, 0\}}=\frac{47 p_1^4 - 470 p_1^2 p_2 +  1445 p_2^2 - 360 p_1 p_3 + 450 p_4}{8},\
\nonumber \\&&
Y_{4,\{1, 0, 1, 1, 0, 0, 0, 0\}}=\frac{43 p_1^4 - 390 p_1^2 p_2 + 873 p_2^2 - 24 p_1 p_3 + 162 p_4}{8},\
\nonumber \\&&
Y_{4,\{0, 0, 0, 2, 0, 0, 0, 0\}}=\frac{19 p_1^4 - 190 p_1^2 p_2 + 553 p_2^2 - 104 p_1 p_3 + 130 p_4}{8},\
\nonumber \\&&
Y_{4,\{3, 0, 0, 0, 1, 0, 0, 0\}}=\frac{43 p_1^4 - 414 p_1^2 p_2 + 1161 p_2^2 - 360 p_1 p_3 + 1242 p_4}{8},\
\nonumber \\&&
Y_{4,\{1, 1, 0, 0, 1, 0, 0, 0\}}=\frac{35 p_1^4 - 286 p_1^2 p_2 + 921 p_2^2 - 488 p_1 p_3 +  690 p_4}{8},\
\nonumber \\&&
Y_{4,\{0, 0, 1, 0, 1, 0, 0, 0\}}=\frac{31 p_1^4 - 222 p_1^2 p_2 + 405 p_2^2 - 24 p_1 p_3 + 138 p_4}{8},\
\nonumber \\&&
Y_{4,\{2, 0, 0, 0, 0, 1, 0, 0\}}=\frac{19 p_1^4 - 142 p_1^2 p_2 +  585 p_2^2 - 584 p_1 p_3 + 1170 p_4}{8},\
\nonumber \\&&
Y_{4,\{0, 1, 0, 0, 0, 1, 0, 0\}}=\frac{3 p_1^4 - 30 p_1^2 p_2 + 465 p_2^2 - 520 p_1 p_3 + 650 p_4}{8},\
\nonumber \\&&
Y_{4,\{1, 0, 0, 0, 0, 0, 1, 0\}}=\frac{-p_1^4 + 34 p_1^2 p_2 + 261 p_2^2 - 600 p_1 p_3 + 858 p_4}{8},\
\nonumber \\&&
Y_{4,\{0, 0, 0, 0, 0, 0, 0, 1\}}=-\frac{3 (3 p_1^4 - 30 p_1^2 p_2 - 39 p_2^2 + 152 p_1 p_3 - 190 p_4)}{8}.
\nonumber 
\end{eqnarray}

\section{Coefficients $C_{n,\bf k}$ for $Y_{n-m}({\bf x}^m)=0$}
\renewcommand{\theequation}{\thesection\arabic{equation}}
\setcounter{equation}{0}
\label{coeff_C}
Here we present the relations for the coefficients $C_{n,\bf k}$ for which the 
homogeneous symmetric polynomials $Y_{n-m}({\bf x}^m)$ vanish for all $ n \le 6$.
\bea
&&
C_{2, \{0, 1\}} = - C_{2, \{2, 0\}}/3,
\nonumber\\
&&
C_{3, \{1, 1, 0\}} = -C_{3, \{3, 0, 0\}}, \
C_{3, \{0, 0, 1\}} = 0,
\nonumber\\
&&
C_{4, \{0, 2, 0, 0\}} = -(5\; C_{4, \{2, 1, 0, 0\}} + 9\; C_{4, \{4, 0, 0, 0\}})/3,
\
C_{4, \{1, 0, 1, 0\}} = -2\; C_{4, \{2, 1, 0, 0\}} - 4\; C_{4, \{4, 0, 0, 0\}},
\nonumber\\
&&
C_{4, \{0, 0, 0, 1\}} = 2(22\; C_{4, \{2, 1, 0, 0\}} + 45\; C_{4, \{4, 0, 0, 0\}})/15,
\nonumber\\
&&
C_{5, \{1, 2, 0, 0, 0\}} = -(9\; C_{5, \{3, 1, 0, 0, 0\}} + 25\; C_{5, \{5, 0, 0, 0, 0\}})/3, \
C_{5, \{2, 0, 1, 0, 0\}} = -3\; C_{5, \{3, 1, 0, 0, 0\}} -  10\; C_{5, \{5, 0, 0, 0, 0\}}, 
\nonumber\\
&&  
C_{5, \{0, 1, 1, 0, 0\}} = 5(3\; C_{5, \{3, 1, 0, 0, 0\}} + 10\; C_{5, \{5, 0, 0, 0, 0\}})/3, \
C_{5, \{0, 0, 0, 0, 1\}} = -2(3\; C_{5, \{3, 1, 0, 0, 0\}} + 10\; C_{5, \{5, 0, 0, 0, 0\}}),
\nonumber\\
&&
C_{5, \{1, 0, 0, 1, 0\}} =  2(9\; C_{5, \{3, 1, 0, 0, 0\}} + 31\; C_{5, \{5, 0, 0, 0, 0\}})/3, 
\nonumber\\
&&
C_{6, \{0, 3, 0, 0, 0, 0\}} =- (27\; C_{6, \{2, 2, 0, 0, 0, 0\}} +115\; C_{6, \{4, 1, 0, 0, 0, 0\}} + 445\; C_{6, \{6, 0, 0, 0, 0, 0\}})/9, 
\nonumber\\
&&
C_{6, \{3, 0, 1, 0, 0, 0\}} = -4 (C_{6, \{4, 1, 0, 0, 0, 0\}} +5\; C_{6, \{6, 0, 0, 0, 0, 0\}}), 
\nonumber\\
&&
C_{6, \{1, 1, 1, 0, 0, 0\}} = -4 (C_{6, \{2, 2, 0, 0, 0, 0\}} + 4\; C_{6, \{4, 1, 0, 0, 0, 0\}}), 
\nonumber\\
&&
C_{6, \{0, 0, 2, 0, 0, 0\}} = 2 (C_{6, \{2, 2, 0, 0, 0, 0\}} + 2\; C_{6, \{4, 1, 0, 0, 0, 0\}} +  5\; C_{6, \{6, 0, 0, 0, 0, 0\}}), 
\nonumber\\
&&
C_{6, \{2, 0, 0, 1, 0, 0\}} = (-7\; C_{6, \{2, 2, 0, 0, 0, 0\}} +  18\; C_{6, \{4, 1, 0, 0, 0, 0\}} + 135\; C_{6, \{6, 0, 0, 0, 0, 0\}})/5,
\nonumber\\
&& 
C_{6, \{0, 1, 0, 1, 0, 0\}} = (207\; C_{6, \{2, 2, 0, 0, 0, 0\}} + 832\; C_{6, \{4, 1, 0, 0, 0, 0\}} +  3115\; C_{6, \{6, 0, 0, 0, 0, 0\}})/15, 
\nonumber\\
&&
C_{6, \{1, 0, 0, 0, 1, 0\}} =  24 (C_{6, \{2, 2, 0, 0, 0, 0\}} + C_{6, \{4, 1, 0, 0, 0, 0\}})/5, 
\nonumber\\
&&
C_{6, \{0, 0, 0, 0, 0, 1\}} = -2(2187\; C_{6, \{2, 2, 0, 0, 0, 0\}} +8102\; C_{6, \{4, 1, 0, 0, 0, 0\}} +29615\; C_{6, \{6, 0, 0, 0, 0, 0\}})/315.
\nonumber
\eea

\end{document}